\documentclass[11pt,reqno]{amsart}
\usepackage{graphicx}
\usepackage{verbatim}
\usepackage{textcomp}
\usepackage{amssymb}
\usepackage{cite}
\usepackage{amsmath}
\usepackage{latexsym}
\usepackage{amscd}
\usepackage{amsthm}
\usepackage{mathrsfs}
\usepackage{xypic}
\usepackage{bm}
\usepackage{url}
\usepackage{hyperref}

\newtheorem{thm}{Theorem}[section]
\newtheorem{corr}[thm]{Corollary}
\newtheorem{lem}[thm]{Lemma}

\theoremstyle{definition}

\theoremstyle{remark}
\newtheorem{rem}{Remark}[section]
\numberwithin{equation}{section}
\begin{document}
\title[Rigidity of complete Riemannian manifolds]
{Rigidity of complete Riemannian manifolds with vanishing Bach tensor}

\author{Bingqing Ma}

\author{Guangyue Huang }
\address{Department of Mathematics, Henan Normal
University, Xinxiang 453007, P.R. China}
\email{bqma@henannu.edu.cn }
\email{hgy@henannu.edu.cn }

\thanks{The research of the first author is supported by NSFC(No. 11401179) and the second author is
supported by NSFC(Nos. 11371018, 11671121).
}

\maketitle

\begin{abstract}

For complete Riemannian manifolds with vanishing Bach tensor and positive constant scalar curvature, we provide a rigidity theorem characterized by some pointwise inequalities. Furthermore, we prove some rigidity results under an inequality involving $L^{\frac{n}{2}}$-norm of the Weyl curvature, the traceless Ricci curvature and the Sobolev constant.

\end{abstract}



\section{Introduction}

In order to study conformal relativity, R. Bach \cite{Bach21} in early
1920s' introduced the following Bach tensor
\begin{equation}\label{1-Sec-1}
B_{ij}=\frac{1}{n-3}W_{ikjl,lk}+\frac{1}{n-2}W_{ikjl}R_{kl},
\end{equation}
where $n\geq4$, $W_{ijkl}$ denotes the Weyl curvature.
A Bach tensor of the metric $g$ is called a vanishing Bach tensor if $B_{ij}=0$. The authors in \cite{Chu2011,Chu2012,Kim2010} consider complete noncompact Riemannian manifolds with vanishing Bach tensor and prove that $M^n$ is of constant curvature if the $L^2$-norm of traceless Riemannian curvature tensor is small. In \cite{Kim2011}, Kim studied complete noncompact Riemannian manifolds with harmonic curvature and positive Sobolev constant, he obtained that $M^n$, $n\geq5$, is Einstein if the $L^2$-norm of the Weyl curvature and the traceless Ricci curvature are small enough.

The aim of this paper is to achieve some rigidity results for complete Riemannian manifolds with vanishing Bach tensor.
In order to state our results, throughout this paper, we always denote by $R$, $\mathring{R}_{ij}$ the scalar curvature and the traceless Ricci curvature of $M^{n}(n\geq4)$, respectively.

\begin{thm}\label{thm1-1}
Let $(M^{n},g)$ be a complete manifold with vanishing Bach tensor,
positive constant scalar curvature and
\begin{equation}\label{1Th-2}
\int_M|\mathring{R}_{ij}|^2<\infty.
\end{equation}
If
\begin{equation}\label{1Th-1}\aligned
\Big|W+\frac{n}{\sqrt{8n}(n-2)}\mathring{{\rm Ric}} \mathbin{\bigcirc\mkern-15mu\wedge} g\Big|\leq\frac{R}{\sqrt{2(n-1)(n-2)}},
\endaligned\end{equation}
then $M^{n}$ is Einstein. In particular, when $n=4,5$, $M^{n}$ is of constant positive sectional curvature.

\end{thm}

Recall that the Sobolev constant $Q_g(M)$ is defined by
\begin{equation}\label{1-Sec-4}
\aligned
Q_g(M)=\inf\limits_{0\neq u\in C_0^{\infty}(M)}\frac{\int_M\Big(|\nabla u|^2+\frac{n-2}{4(n-1)}R_gu^2\Big)}{(\int_M|u|^{\frac{2n}{n-2}})^{\frac{n-2}{n}}}.
\endaligned\end{equation}
Moreover, there exist complete noncompact manifolds with negative scalar curvature which have positive Sobolev constant. For example, any simply connected complete manifold with $W_{ijkl}=0$ has positive Sobolev constant (see \cite{Schoen-yau1988}).
Moreover, it is easy to see from \eqref{1-Sec-4}, for any $u$,
\begin{equation}\label{1-Sec-5}
\aligned
Q_g(M)\Big(\int_M|u|^{\frac{2n}{n-2}}\Big)^{\frac{n-2}{n}}\leq\int_M\Big(|\nabla u|^2+\frac{n-2}{4(n-1)}R_gu^2\Big).
\endaligned\end{equation}

With the help of \eqref{1-Sec-5}, we can achieve the following rigidity result:

\begin{thm}\label{thm1-2}
Let $(M^{n},g)$ be a complete manifold with vanishing Bach tensor,
constant scalar curvature, $Q_g(M)>0$ and
\begin{equation}\label{2Th-2}\aligned
\int_M|\mathring{R}_{ij}|^2<\infty.
\endaligned\end{equation}

(1) If $n\geq7$ and
\begin{equation}\label{2Th-1}\aligned
\Big(\int_M\Big|W+\frac{\sqrt{n}}{\sqrt{8}(n-2)}\mathring{{\rm Ric}} \mathbin{\bigcirc\mkern-15mu\wedge} g\Big|^{\frac{n}{2}}\Big)^{\frac{2}{n}}<\frac{2}{n-2}\sqrt{\frac{2(n-1)}{n-2}}Q_g(M),
\endaligned\end{equation}
then $M^{n}$ is Einstein;

(2) If $4\leq n\leq 6$, $R\geq0$ and
\begin{equation}\label{2Th-3}\aligned
\Big(\int_M\Big|W+\frac{\sqrt{n}}{\sqrt{8}(n-2)}\mathring{{\rm Ric}} \mathbin{\bigcirc\mkern-15mu\wedge} g\Big|^{\frac{n}{2}}\Big)^{\frac{2}{n}}<\sqrt{\frac{n-1}{2(n-2)}}Q_g(M),
\endaligned\end{equation}
then $M^{n}$ is Einstein. In particular, for $M^{n}(n=4,5)$ with positive constant scalar curvature, it must be of constant sectional curvature.
\end{thm}

It is well-known that there is no complete noncompact Einstein manifold with positive constant scalar curvature. Hence, the following results follow easily:

\begin{corr}\label{1corr1-1}
Suppose that $(M^n,g)$ is a complete noncompact Riemannian manifold with vanishing Bach tensor and
positive constant scalar curvature. If \eqref{1Th-1} holds, then we have
\begin{align}\label{corr-Int-1}
\int_M|\mathring{\rm Ric}|^2=\infty.
\end{align}

\end{corr}

\begin{corr}\label{2corr1-1}
Let $(M^{n},g)$ be a complete noncompact manifold with vanishing Bach tensor and positive constant scalar curvature.
If either $n\geq7$ and \eqref{2Th-1}, or $4\leq n\leq 6$ and \eqref{2Th-3} holds, then we have
\begin{equation}\label{1corr-1}
\int_M|\mathring{R}_{ij}|^2=\infty.
\end{equation}

\end{corr}

\begin{rem}
The authors in \cite{Yuan2017,Huang2017} obtained some rigidity results
similar to our theorems for compact manifolds with vanishing Bach tensor. Our theorems can be seen as a generalization to complete manifolds.
\end{rem}

\begin{rem}
In \cite[Theorem 2]{Kim2010}, Kim proved, for complete noncompact manifold $M^4$ with vanishing Bach tensor and nonnegative constant scalar curvature, that if there exists a constant $c_0$ small enough such that
\begin{equation}\label{rem-1}
 \int_M(|W|^2+|\mathring{R}_{ij}|^2)\leq c_0,
\end{equation}
then $M^4$ is Einstein. When $n=4$, our formula \eqref{2Th-3} becomes
\begin{equation}\label{rem-2}\aligned
\int_M\Big|W+\frac{1}{2\sqrt{2}}\mathring{{\rm Ric}} \mathbin{\bigcirc\mkern-15mu\wedge} g\Big|^2<\frac{3}{4}Q_g^2(M),
\endaligned\end{equation}
which, compared with \eqref{rem-1}, shows that our Theorem \ref{thm1-2} gives an upper bound of $c_0$, in some sense.
Moreover, we also prove that for the upper bound $\frac{3}{4}Q_g^2(M)$ given by \eqref{rem-2}, $M^4$ must be of constant sectional curvature if it has positive constant scalar curvature.

\end{rem}

\section{Some lemmas}

Recall that the Weyl curvature $W_{ijkl}$ is relate to the Riemannian curvature $R_{ijkl}$ by
\begin{equation}\label{2-Sec-1}\aligned
W_{ijkl}=&R_{ijkl}-\frac{1}{n-2}(R_{ik}g_{jl}-R_{il}g_{jk}
+R_{jl}g_{ik}-R_{jk}g_{il}).
\endaligned
\end{equation}
By virtue of $\mathring{R}_{ij}=R_{ij}-\frac{R}{n}g_{ij}$, \eqref{2-Sec-1} can be written as
\begin{equation}\label{add2-Sec-1}\aligned
W_{ijkl}=&R_{ijkl}-\frac{1}{n-2}(\mathring{R}_{ik}g_{jl}-\mathring{R}_{il}g_{jk}
+\mathring{R}_{jl}g_{ik}-\mathring{R}_{jk}g_{il})\\
&-\frac{R}{n(n-1)}(g_{ik}g_{jl}-g_{il}g_{jk}).
\endaligned
\end{equation}
Since the divergence of the Weyl curvature tensor is related to the Cotton tensor by
\begin{equation}\label{2-Sec-2}
W_{ijkl,l}=-\frac{n-3}{n-2}C_{ijk},
\end{equation}
where the Cotton tensor is given by
\begin{equation}\label{2-Sec-3}\aligned
C_{ijk}=&R_{kj,i}-R_{ki,j}-\frac{1}{2(n-1)}(R_{,i}g_{jk}-R_{,j}g_{ik})\\
=&\mathring{R}_{kj,i}-\mathring{R}_{ki,j}+\frac{n-2}{2n(n-1)}(R_{,i}g_{jk}-R_{,j}g_{ik}),
\endaligned
\end{equation}
the formula \eqref{1-Sec-1} reduces to
\begin{equation}\label{2-Sec-4}
B_{ij}=\frac{1}{n-2}(C_{kij,k}+W_{ikjl}R^{kl}).
\end{equation}

We denote by $p\in M^n$ and $B_r$ a fixed point and the geodesic sphere of $M^n$ of radius $r$ centered at $p$, respectively.  Let $\phi_r$ be the nonnegative cut-off function defined on $M^n$ satisfying
\begin{equation}\label{cutoff}
\phi_r=
\begin{cases}1,\ \ \quad {\rm on}\ B_r\\
0,\ \ \quad {\rm on}\ M^n\backslash B_{r+1}
\end{cases}
\end{equation}
with $|\nabla \phi_r|\leq2$ on $B_{r+1}\backslash B_{r}$.

Inspired by \cite[Lemma 2.2]{Yuan2017}, we give the following estimate with respect to $\mathring{R}_{ij}$ on a complete noncompact Riemannian manifold:

\begin{lem}\label{2lem-1}
Let $(M^n,g)$ be a complete noncompact Riemannian manifold with constant scalar curvature. Then for any $\theta\in \mathbb{R}$, we have
\begin{equation}\label{2-Sec-5}\aligned
&\int_M|\nabla \mathring{R}_{ij}|^2\phi_{r}^2\\
\geq&\frac{2\theta}{\theta^2+1+\epsilon_1}\int_M\Big(W_{ijkl}\mathring{R}_{ik}\mathring{R}_{jl}-\frac{n}{n-2}\mathring{R}_{ij}\mathring{R}_{jk}
\mathring{R}_{ki}-\frac{1}{n-1}R|\mathring{R}_{ij}|^2\Big)\phi_{r}^2\\
&-\frac{4\theta^2}{\epsilon_1(\theta^2+1+\epsilon_1)}\int_M|\mathring{R}_{ij}|^2|\nabla \phi_{r}|^2,
\endaligned\end{equation}
where $\epsilon_1$ is a positive constant.
\end{lem}

\proof  By a direct calculation, we have
\begin{equation}\aligned\label{2-Sec-6}
0\leq&\int_M|\mathring{R}_{kj,i}-\theta\mathring{R}_{ki,j}|^2\phi_{r}^2\\
=&(\theta^2+1)\int_M|\nabla \mathring{R}_{ij}|^2\phi_{r}^2-2\theta\int_M\mathring{R}_{kj,i}\mathring{R}_{ki,j}\phi_{r}^2.
\endaligned\end{equation}
Using the Ricci identity and \eqref{add2-Sec-1}, we have
\begin{equation}\aligned\label{2-Sec-7}
\mathring{R}_{kj,ij}\mathring{R}_{ki}=&(\mathring{R}_{kj,ji}+\mathring{R}_{pj}R_{pkij}+\mathring{R}_{kp}R_{pjij})\mathring{R}_{ki}\\
=&\mathring{R}_{pj}\mathring{R}_{ki}R_{pkij}
+\mathring{R}_{kp}\mathring{R}_{ki}\mathring{R}_{pi}+\frac{R}{n}|\mathring{R}_{ij}|^2\\
=&-W_{ijkl}\mathring{R}_{ik}\mathring{R}_{jl}+\frac{n}{n-2}\mathring{R}_{ij}\mathring{R}_{jk}
\mathring{R}_{ki}+\frac{1}{n-1}R|\mathring{R}_{ij}|^2,
\endaligned\end{equation}
where we used the second Bianchi identity and $\mathring{R}_{ij,j}=\frac{n-2}{2n}R_{,i}=0$ from the fact that $R$ is constant.
Hence, we obtain
\begin{equation}\aligned\label{2-Sec-8}
-2\theta\int_M\mathring{R}_{kj,i}\mathring{R}_{ki,j}\phi_{r}^2
=&2\theta\int_M\mathring{R}_{kj,ij}\mathring{R}_{ki}\phi_{r}^2
+2\theta\int_M\mathring{R}_{kj,i}\mathring{R}_{ki}(\phi_{r}^2)_{j}\\
=&-2\theta\int_M\Big(W_{ijkl}\mathring{R}_{ik}\mathring{R}_{jl}-\frac{n}{n-2}\mathring{R}_{ij}\mathring{R}_{jk}
\mathring{R}_{ki}-\frac{1}{n-1}R|\mathring{R}_{ij}|^2\Big)\phi_{r}^2\\
&+2\theta\int_M\mathring{R}_{kj,i}\mathring{R}_{ki}(\phi_{r}^2)_{j}\\
\leq&-2\theta\int_M\Big(W_{ijkl}\mathring{R}_{ik}\mathring{R}_{jl}-\frac{n}{n-2}\mathring{R}_{ij}\mathring{R}_{jk}
\mathring{R}_{ki}-\frac{1}{n-1}R|\mathring{R}_{ij}|^2\Big)\phi_{r}^2\\
&+\epsilon_1\int_M|\nabla \mathring{R}_{ij}|^2\phi_{r}^2+\frac{4\theta^2}{\epsilon_1}\int_M|\mathring{R}_{ij}|^2|\nabla \phi_{r}|^2.
\endaligned\end{equation}
Applying \eqref{2-Sec-8} into \eqref{2-Sec-6} yields the desired estimate \eqref{2-Sec-5}.
\endproof

\begin{lem}\label{2lem-2}
Let $(M^n,g)$ be a complete noncompact Riemannian manifold with constant scalar curvature. If the Bach tensor is flat, we have
\begin{equation}\label{2-Sec-9}\aligned
\int_M|\nabla \mathring{R}_{ij}|^2\phi_{r}^2
\leq&\frac{1}{1-\epsilon_2}\int_M\Big(2W_{ijkl}\mathring{R}_{jl}\mathring{R}_{ik}
-\frac{n}{n-2}\mathring{R}_{ij}\mathring{R}_{jk}
\mathring{R}_{ki}\\
&-\frac{1}{n-1}R|\mathring{R}_{ij}|^2\Big)\phi_{r}^2+\frac{1}{\epsilon_2(1-\epsilon_2)}\int_M|\mathring{R}_{ij}|^2|\nabla \phi_{r}|^2,
\endaligned\end{equation}
where $\epsilon_2\in(0,1)$ is a constant.
\end{lem}

\proof Using the formula \eqref{add2-Sec-1}, we can derive
\begin{equation}\label{2-Sec-10}\aligned
\mathring{R}_{kl}R_{ikjl}=&\mathring{R}_{kl}W_{ikjl}+\frac{1}{n-2}(|\mathring{R}_{ij}|^2g_{ij}
-2\mathring{R}_{ik}\mathring{R}_{jk})-\frac{1}{n(n-1)}R \mathring{R}_{ij},
\endaligned\end{equation}
which shows
\begin{equation}\label{2-Sec-11}\aligned
\mathring{R}_{kj,ik}=&\mathring{R}_{kj,ki}+\mathring{R}_{lj}R_{lkik}+\mathring{R}_{kl}R_{ljik}\\
=&\mathring{R}_{ik}\mathring{R}_{jk}+\frac{1}{n}R\mathring{R}_{ij}-\Big[\mathring{R}_{kl}W_{ikjl}
\\
&+\frac{1}{n-2}(|\mathring{R}_{ij}|^2g_{ij}-2\mathring{R}_{ik}\mathring{R}_{jk})
-\frac{1}{n(n-1)}R \mathring{R}_{ij}\Big]\\
=&\frac{n}{n-2}\mathring{R}_{ik}\mathring{R}_{jk}+\frac{1}{n-1}R\mathring{R}_{ij}
-\mathring{R}_{kl}W_{ikjl}-\frac{1}{n-2}|\mathring{R}_{ij}|^2g_{ij}.
\endaligned\end{equation}
Thus, from \eqref{2-Sec-3} and \eqref{2-Sec-11}, we have
\begin{equation}\label{2-Sec-12}\aligned
C_{kij,k}=&\Delta \mathring{R}_{ij}-\mathring{R}_{kj,ik}\\
=&\Delta \mathring{R}_{ij}-\Big(\frac{n}{n-2}\mathring{R}_{ik}\mathring{R}_{jk}
+\frac{1}{n-1}R\mathring{R}_{ij}-\mathring{R}_{kl}W_{ikjl}\\
&-\frac{1}{n-2}|\mathring{R}_{ij}|^2g_{ij}\Big)
\endaligned\end{equation}
and
\begin{equation}\label{2-Sec-13}\aligned
0=&(n-2)B_{ij}\mathring{R}_{ij}\\
=&C_{kij,k}\mathring{R}_{ij}+W_{ikjl}\mathring{R}_{ij}\mathring{R}_{kl}\\
=&\mathring{R}_{ij}\Delta \mathring{R}_{ij}-\frac{n}{n-2}\mathring{R}_{ij}\mathring{R}_{jk}\mathring{R}_{ki}
-\frac{1}{n-1}R|\mathring{R}_{ij}|^2+2W_{ikjl}\mathring{R}_{ij}\mathring{R}_{kl},
\endaligned\end{equation}
which gives
\begin{equation}\label{add2-Sec-14}
\mathring{R}_{ij}\Delta \mathring{R}_{ij}=\frac{n}{n-2}\mathring{R}_{ij}\mathring{R}_{jk}\mathring{R}_{ki}
+\frac{1}{n-1}R|\mathring{R}_{ij}|^2-2W_{ikjl}\mathring{R}_{ij}\mathring{R}_{kl}.
\end{equation}
Thus,
\begin{equation}\label{2-Sec-14}\aligned
\int_M|\nabla \mathring{R}_{ij}|^2\phi_{r}^2=&-\int_M\mathring{R}_{ij}\Delta\mathring{R}_{ij}\phi_{r}^2
-\int_M\mathring{R}_{ij}\mathring{R}_{ij,k}(\phi_{r}^2)_{k}\\
=&\int_M\Big(2W_{ijkl}\mathring{R}_{jl}\mathring{R}_{ik}
-\frac{n}{n-2}\mathring{R}_{ij}\mathring{R}_{jk}
\mathring{R}_{ki}\\
&-\frac{1}{n-1}R|\mathring{R}_{ij}|^2\Big)\phi_{r}^2-\int_M\mathring{R}_{ij}\mathring{R}_{ij,k}(\phi_{r}^2)_{k}\\
\leq&\int_M\Big(2W_{ijkl}\mathring{R}_{jl}\mathring{R}_{ik}
-\frac{n}{n-2}\mathring{R}_{ij}\mathring{R}_{jk}
\mathring{R}_{ki}\\
&-\frac{1}{n-1}R|\mathring{R}_{ij}|^2\Big)\phi_{r}^2+\epsilon_2\int_M|\nabla \mathring{R}_{ij}|^2\phi_{r}^2+\frac{1}{\epsilon_2}\int_M|\mathring{R}_{ij}|^2|\nabla \phi_{r}|^2.
\endaligned\end{equation}
We complete the proof of Lemma \ref{2lem-2}.\endproof

The following two lemmas come from \cite{Huang2017} (for more details, see \cite{Huang2017jmaa,Catino2016,FuXiao2015,MHLC2018,MH2017}. For a proof of Lemma \ref{2lem-4}, we refer to \cite{HCL2018}):

\begin{lem}\label{2lem-3}
On every Riemannian manifold $(M^n,g)$, for any $\lambda\in \mathbb{R}$, the following estimate holds
\begin{equation}\label{2-Sec-18}\aligned
\Big|-&W_{ijkl}\mathring{R}_{jl}\mathring{R}_{ik}
+\lambda\mathring{R}_{ij}\mathring{R}_{jk}\mathring{R}_{ki}\Big|\\
\leq&
\sqrt{\frac{n-2}{2(n-1)}}\Big(|W|^2+\frac{2(n-2)\lambda^2}{n}|\mathring{R}_{ij}|^2
\Big)^{\frac{1}{2}}|\mathring{R}_{ij}|^2\\
=&\sqrt{\frac{n-2}{2(n-1)}}\Big|W+\frac{\lambda}{\sqrt{2n}}\mathring{{\rm Ric}} \mathbin{\bigcirc\mkern-15mu\wedge} g\Big||\mathring{R}_{ij}|^2.
\endaligned\end{equation}

\end{lem}

\begin{lem}\label{2lem-4}
On every Einstein manifold $(M^n,g)$, we have
\begin{align}\label{2-Sec-19}
\frac{1}{2}\Delta |W|^2\geq&\frac{n+1}{n-1}|\nabla |W||^2+\frac{2}{n}R|W|^2-2C_n |W|^3,
\end{align}
where $C_n$ is defined by
\begin{equation}\label{2-Sec-20}
C_n=\begin{cases}\frac{\sqrt{6}}{4},&\text{if\ }n=4;\\
\frac{4\sqrt{10}}{15},&\text{if\ }n=5;\\
\frac{n-2}{\sqrt{n(n-1)}} +\frac{n^2-n-4}{2\sqrt{(n-2)(n-1)n(n+1)}},&\text{if\ }n\geq 6.
\end{cases}
\end{equation}
In particular, if the scalar curvature of Einstein metric $g$ is positive,
then it is of constant positive sectional curvature, provided either
\begin{align}\label{Integ-1}
C_n|W|<\frac{1}{n}R,
\end{align}
or

(1) for $n\neq5$,
\begin{equation}\label{2-Pinchconstant}
\Big(\int_M|W|^{\frac{n}{2}}\Big)^{\frac{2}{n}}<E_n\,Q_g(M),
\end{equation}
where $E_n$ is given by
\begin{align}\label{En}
E_n=&\begin{cases}\sqrt{6},&\text{if\ }n=4;\\
\frac{4(n-1)}{n(n-2)}\Big(\frac{n-2}{\sqrt{n(n-1)}} +\frac{n^2-n-4}{2\sqrt{(n-2)(n-1)n(n+1)}}\Big)^{-1},&\text{if\ }n\geq 6;
\end{cases}
\end{align}

(2) for $n=5$,
\begin{equation}\label{2-Pinchconstant-2}
\Big(\int_M|W|^{\frac{5}{2}}\Big)^{\frac{2}{5}}\leq\frac{2\sqrt{15}-4}{\sqrt{10}}\,Q_g(M).
\end{equation}

\end{lem}

\section{Proof of theorems}

\subsection{Proof of theorem \ref{thm1-1}}
By combining the estimates \eqref{2-Sec-5} with \eqref{2-Sec-9}, we derive
\begin{equation}\label{3-Sec-1}\aligned
&\frac{[4\epsilon_2(1-\epsilon_2)+\epsilon_1]\theta^2+\epsilon_1(1+\epsilon_1)}{\epsilon_1\epsilon_2
(1-\epsilon_2)(\theta^2+1+\epsilon_1)}\int_M|\mathring{R}_{ij}|^2|\nabla \phi_{r}|^2\\
\geq&-2\frac{\theta^2-(1-\epsilon_2)\theta+(1+\epsilon_1)}{(1-\epsilon_2)(\theta^2+1+\epsilon_1)}
\int_MW_{ijkl}\mathring{R}_{jl}\mathring{R}_{ik}\phi_{r}^2\\
&+\frac{\theta^2-2(1-\epsilon_2)\theta+(1+\epsilon_1)}{(1-\epsilon_2)(\theta^2+1+\epsilon_1)}\frac{n}{n-2}
\int_M\mathring{R}_{ij}\mathring{R}_{jk}\mathring{R}_{ki}\phi_{r}^2\\
&+\frac{\theta^2-2(1-\epsilon_2)\theta+(1+\epsilon_1)}{(1-\epsilon_2)(\theta^2+1+\epsilon_1)}\frac{1}{n-1}
\int_M R|\mathring{R}_{ij}|^2\phi_{r}^2.
\endaligned\end{equation}
For all $\epsilon_2\in(0,1)$, we have \begin{equation}\label{3-Sec-2}\theta^2-2(1-\epsilon_2)\theta+(1+\epsilon_1)=[\theta-(1-\epsilon_2)]^2+(1+\epsilon_1)
-(1-\epsilon_2)^2>0,\end{equation}
and hence \eqref{3-Sec-1} reduces to
\begin{equation}\label{3-Sec-3}\aligned
&\frac{[4\epsilon_2(1-\epsilon_2)+\epsilon_1]\theta^2+\epsilon_1(1+\epsilon_1)}{\epsilon_1\epsilon_2
[\theta^2-2(1-\epsilon_2)\theta+(1+\epsilon_1)]
}\int_M|\mathring{R}_{ij}|^2|\nabla \phi_{r}|^2\\
\geq&-\frac{2[\theta^2-(1-\epsilon_2)\theta+(1+\epsilon_1)]}{\theta^2-2(1-\epsilon_2)\theta+(1+\epsilon_1)}
\int_MW_{ijkl}\mathring{R}_{jl}\mathring{R}_{ik}\phi_{r}^2\\
&+\frac{n}{n-2}\int_M\mathring{R}_{ij}\mathring{R}_{jk}\mathring{R}_{ki}\phi_{r}^2+\frac{1}{n-1}\int_M R|\mathring{R}_{ij}|^2\phi_{r}^2.
\endaligned\end{equation}
Using \eqref{2-Sec-18}, we have
\begin{equation}\label{3-Sec-4}\aligned
-&\frac{2[\theta^2-(1-\epsilon_2)\theta+(1+\epsilon_1)]}{\theta^2-2(1-\epsilon_2)\theta+(1+\epsilon_1)}
W_{ijkl}\mathring{R}_{jl}\mathring{R}_{ik}+\frac{n}{n-2}\mathring{R}_{ij}\mathring{R}_{jk}\mathring{R}_{ki}\\
\geq&-\sqrt{\frac{n-2}{2(n-1)}}\Big|\frac{2[\theta^2-(1-\epsilon_2)\theta+(1+\epsilon_1)]}{\theta^2-2(1-\epsilon_2)\theta
+(1+\epsilon_1)}W
+\frac{n}{\sqrt{2n}(n-2)}\mathring{{\rm Ric}} \mathbin{\bigcirc\mkern-15mu\wedge} g\Big||\mathring{R}_{ij}|^2.
\endaligned\end{equation}
Applying \eqref{3-Sec-4} into \eqref{3-Sec-3} gives
\begin{equation}\label{3-Sec-5}\aligned
&\frac{[4\epsilon_2(1-\epsilon_2)+\epsilon_1]\theta^2+\epsilon_1(1+\epsilon_1)}{\epsilon_1\epsilon_2
[\theta^2-2(1-\epsilon_2)\theta+(1+\epsilon_1)]
}\int_M|\mathring{R}_{ij}|^2|\nabla \phi_{r}|^2\\
\geq&\int_M\Bigg[-\sqrt{\frac{n-2}{2(n-1)}}\Big|\frac{2[\theta^2-(1-\epsilon_2)\theta+(1+\epsilon_1)]}{\theta^2-2(1-\epsilon_2)\theta
+(1+\epsilon_1)}W\\
&+\frac{n}{\sqrt{2n}(n-2)}\mathring{{\rm Ric}} \mathbin{\bigcirc\mkern-15mu\wedge} g\Big|+\frac{1}{n-1} R\Bigg]|\mathring{R}_{ij}|^2\phi_{r}^2.
\endaligned\end{equation}
Now, we fixed $\epsilon_1$ and $\epsilon_2$ and minimize the coefficient of $W$ with respect to the function $\theta$ by taking
\begin{equation}\label{3-Sec-6}
\theta=-\sqrt{1+\epsilon_1},
\end{equation}
then \eqref{3-Sec-5} becomes
\begin{equation}\label{3-Sec-7}\aligned
&\frac{[2\epsilon_2(1-\epsilon_2)+\epsilon_1]\sqrt{1+\epsilon_1}}{\epsilon_1\epsilon_2
[\sqrt{1+\epsilon_1}+(1-\epsilon_2)]
}\int_M|\mathring{R}_{ij}|^2|\nabla \phi_{r}|^2\\
\geq&\int_M\Bigg[-\sqrt{\frac{n-2}{2(n-1)}}\Big|\Big(1+\frac{1}{1+\frac{1-\epsilon_2}{\sqrt{1+\epsilon_1}}}\Big)W\\
&+\frac{n}{\sqrt{2n}(n-2)}\mathring{{\rm Ric}} \mathbin{\bigcirc\mkern-15mu\wedge} g\Big|+\frac{1}{n-1} R\Bigg]|\mathring{R}_{ij}|^2\phi_{r}^2.
\endaligned\end{equation}
Since $W$ is perpendicular to $\mathring{{\rm Ric}} \mathbin{\bigcirc\mkern-15mu\wedge} g$, for any given $\check{\epsilon}_1,\check{\epsilon}_2$, we have
\begin{equation}\label{3-Sec-8}\aligned
\Big|&\Big(1+\frac{1}{1+\frac{1-\check{\epsilon}_2}{\sqrt{1+\check{\epsilon}_1}}}\Big)W+\frac{n}{\sqrt{2n}(n-2)}\mathring{{\rm Ric}} \mathbin{\bigcirc\mkern-15mu\wedge} g\Big|^2\\
&=\Big(1+\frac{1}{1+\frac{1-\check{\epsilon}_2}{\sqrt{1+\check{\epsilon}_1}}}\Big)^2|W|^2+\frac{n}{2(n-2)^2}|\mathring{{\rm Ric}} \mathbin{\bigcirc\mkern-15mu\wedge} g|^2\\
&<4|W|^2+\frac{n}{2(n-2)^2}|\mathring{{\rm Ric}} \mathbin{\bigcirc\mkern-15mu\wedge} g|^2\\
&=\Big|2W+\frac{n}{\sqrt{2n}(n-2)}\mathring{{\rm Ric}} \mathbin{\bigcirc\mkern-15mu\wedge} g\Big|^2,
\endaligned\end{equation}
which shows that
\begin{equation}\label{3-Sec-9}\aligned
\Big|&\Big(1+\frac{1}{1+\frac{1-\check{\epsilon}_2}{\sqrt{1+\check{\epsilon}_1}}}\Big)W+\frac{n}{\sqrt{2n}(n-2)}\mathring{{\rm Ric}} \mathbin{\bigcirc\mkern-15mu\wedge} g\Big|\\
&<2\Big|W+\frac{n}{\sqrt{8n}(n-2)}\mathring{{\rm Ric}} \mathbin{\bigcirc\mkern-15mu\wedge} g\Big|.
\endaligned\end{equation}
Therefore, under the condition \eqref{1Th-1}, the estimate \eqref{3-Sec-7} gives
\begin{equation}\label{3-Sec-10}\aligned
0\leq&\int_M\Bigg[-\sqrt{\frac{n-2}{2(n-1)}}\Big|\Big(1+\frac{1}{1+\frac{1-\check{\epsilon}_2}{\sqrt{1+\check{\epsilon}_1}}}\Big)W\\
&+\frac{n}{\sqrt{2n}(n-2)}\mathring{{\rm Ric}} \mathbin{\bigcirc\mkern-15mu\wedge} g\Big|+\frac{1}{n-1} R\Bigg]|\mathring{R}_{ij}|^2\phi_{r}^2\\
\leq&\frac{[2\check{\epsilon}_2(1-\check{\epsilon}_2)+\check{\epsilon}_1]\sqrt{1+\check{\epsilon}_1}}{\check{\epsilon}_1\check{\epsilon}_2
[\sqrt{1+\check{\epsilon}_1}+(1-\check{\epsilon}_2)]
}\int_M|\mathring{R}_{ij}|^2|\nabla \phi_{r}|^2.
\endaligned\end{equation}

Since \begin{equation}\label{3-Sec-11}
\int_M|\mathring{R}_{ij}|^2<\infty,
\end{equation}
then we have
\begin{equation}\label{3-Sec-12}\aligned
\int_M|\mathring{R}_{ij}|^2|\nabla \phi_{r}|^2\rightarrow0,
\endaligned\end{equation}
as $r\rightarrow \infty$.
This, together with \eqref{3-Sec-10}, shows that $M^n$ is Einstein. In this case,
\eqref{1Th-1} becomes
\begin{equation}\label{3-Sec-13}\aligned
|W|\leq\frac{R}{\sqrt{2(n-1)(n-2)}},
\endaligned\end{equation}
which yields
\begin{equation}\label{3-Sec-14}\aligned
C_n|W|\leq\frac{C_n}{\sqrt{2(n-1)(n-2)}}R.
\endaligned\end{equation}
It is easy to check that for $n=4,5$, we have
\begin{equation}\label{3-Sec-15}\aligned
\frac{C_n}{\sqrt{2(n-1)(n-2)}}R<\frac{1}{n}R,
\endaligned\end{equation}
which combining with \eqref{Integ-1} shows that $M^n$ is of constant positive sectional curvature.

\subsection{Proof of theorem \ref{thm1-2}}
From \eqref{2-Sec-18}, it is easy to see
\begin{equation}\label{3-Sec-13}\aligned
2W_{ijkl}&\mathring{R}_{jl}\mathring{R}_{ik}-\frac{n}{n-2}
\mathring{R}_{ij}\mathring{R}_{jk}\mathring{R}_{ki}\\
\leq&\sqrt{\frac{2(n-2)}{n-1}}\Big|W+\frac{\sqrt{n}}{\sqrt{8}(n-2)}\mathring{{\rm Ric}} \mathbin{\bigcirc\mkern-15mu\wedge} g\Big||\mathring{R}_{ij}|^2.
\endaligned\end{equation}
Applying \eqref{3-Sec-13} into \eqref{2-Sec-9} and using the Kato inequality, we obtain
\begin{equation}\label{3-Sec-14}
\aligned
\int_M|\nabla |\mathring{R}_{ij}||^2\phi_{r}^2
\leq&\int_M|\nabla \mathring{R}_{ij}|^2\phi_{r}^2\\
\leq&\frac{1}{1-\epsilon_2}\int_M\Bigg[\sqrt{\frac{2(n-2)}{n-1}}\Big|W+\frac{\sqrt{n}}{\sqrt{8}(n-2)}\mathring{{\rm Ric}} \mathbin{\bigcirc\mkern-15mu\wedge} g\Big|\\
&-\frac{1}{n-1}R\Bigg]|\mathring{R}_{ij}|^2\phi_{r}^2+\frac{1}{\epsilon_2(1-\epsilon_2)}\int_M|\mathring{R}_{ij}|^2|\nabla \phi_{r}|^2.
\endaligned\end{equation}
Taking $u=|\mathring{R}_{ij}|\phi_{r}$ in \eqref{1-Sec-5} and applying \eqref{3-Sec-14} yield
\begin{equation}\label{3-Sec-15}
\aligned
Q_g(M)&\Big(\int_M(|\mathring{R}_{ij}|\phi_{r})^{\frac{2n}{n-2}}\Big)^{\frac{n-2}{n}}\\
\leq&\int_M\Big(|\nabla (|\mathring{R}_{ij}|\phi_{r})|^2+\frac{n-2}{4(n-1)}R|\mathring{R}_{ij}|^2\phi_{r}^2\Big)\\
\leq&(1+\epsilon_3)\int_M|\nabla |\mathring{R}_{ij}||^2\phi_{r}^2+\Big(1+\frac{1}{\epsilon_3}\Big)\int_M|\mathring{R}_{ij}|^2|\nabla \phi_{r}|^2\\
&+\frac{n-2}{4(n-1)}\int_MR|\mathring{R}_{ij}|^2\phi_{r}^2\\
\leq&\frac{1+\epsilon_3}{1-\epsilon_2}\sqrt{\frac{2(n-2)}{n-1}}\int_M\Big|W+\frac{\sqrt{n}}{\sqrt{8}(n-2)}\mathring{{\rm Ric}} \mathbin{\bigcirc\mkern-15mu\wedge} g\Big||\mathring{R}_{ij}|^2\phi_{r}^2\\
&+\frac{1}{n-1}\Big[\frac{n-2}{4}-\frac{1+\epsilon_3}{1-\epsilon_2}\Big]\int_MR|\mathring{R}_{ij}|^2\phi_{r}^2\\
&+(1+\epsilon_3)\Big[\frac{1}{\epsilon_3}+\frac{1}{\epsilon_2(1-\epsilon_2)}\Big]\int_M|\mathring{R}_{ij}|^2|\nabla \phi_{r}|^2.
\endaligned\end{equation}
Inserting the following H\"{o}lder inequality
$$\aligned
\int_M&\Big|W+\frac{\sqrt{n}}{\sqrt{8}(n-2)}\mathring{{\rm Ric}} \mathbin{\bigcirc\mkern-15mu\wedge} g\Big||\mathring{R}_{ij}|^2\phi_{r}^2\\
\leq&\Big(\int_M\Big|W+\frac{\sqrt{n}}{\sqrt{8}(n-2)}\mathring{{\rm Ric}} \mathbin{\bigcirc\mkern-15mu\wedge} g\Big|^{\frac{n}{2}}\Big)^{\frac{2}{n}}\Big(\int_M(|\mathring{R}_{ij}|\phi_{r})^{\frac{2n}{n-2}}\Big)^{\frac{n-2}{n}}
\endaligned$$
into \eqref{3-Sec-15} yields
\begin{equation}\label{3-Sec-16}
\aligned
\Bigg[Q_g(M)&-\frac{1+\epsilon_3}{1-\epsilon_2}\sqrt{\frac{2(n-2)}{n-1}}\Big(\int_M\Big|W\\
&+\frac{\sqrt{n}}{\sqrt{8}(n-2)}\mathring{{\rm Ric}} \mathbin{\bigcirc\mkern-15mu\wedge} g\Big|^{\frac{n}{2}}\Big)^{\frac{2}{n}}\Bigg]\Big(\int_M(|\mathring{R}_{ij}|\phi_{r})^{\frac{2n}{n-2}}\Big)^{\frac{n-2}{n}}\\
\leq&\frac{1}{n-1}\Big[\frac{n-2}{4}-\frac{1+\epsilon_3}{1-\epsilon_2}\Big]\int_MR|\mathring{R}_{ij}|^2\phi_{r}^2\\
&+(1+\epsilon_3)\Big[\frac{1}{\epsilon_3}+\frac{1}{\epsilon_2(1-\epsilon_2)}\Big]\int_M|\mathring{R}_{ij}|^2|\nabla \phi_{r}|^2.
\endaligned\end{equation}

Now, we consider the following two cases:

{\it Case one:} When $n\geq7$,  there exist $\tilde{\epsilon}_2,\tilde{\epsilon}_3$ depending only on the dimension $n$ such that
\begin{equation}\label{3-Sec-17}\aligned
\frac{1+\tilde{\epsilon}_3}{1-\tilde{\epsilon}_2}=\frac{n-2}{4}.
\endaligned\end{equation}
In this case, \eqref{3-Sec-16} becomes
\begin{equation}\label{3-Sec-18}\aligned
\Bigg[Q_g(M)&-\frac{n-2}{4}\sqrt{\frac{2(n-2)}{n-1}}\Big(\int_M\Big|W\\
&+\frac{\sqrt{n}}{\sqrt{8}(n-2)}\mathring{{\rm Ric}} \mathbin{\bigcirc\mkern-15mu\wedge} g\Big|^{\frac{n}{2}}\Big)^{\frac{2}{n}}\Bigg]\Big(\int_M(|\mathring{R}_{ij}|\phi_{r})^{\frac{2n}{n-2}}\Big)^{\frac{n-2}{n}}\\
\leq&(1+\tilde{\epsilon}_3)\Big[\frac{1}{\tilde{\epsilon}_3}+\frac{1}{\tilde{\epsilon}_2(1-\tilde{\epsilon}_2)}\Big]\int_M|\mathring{R}_{ij}|^2|\nabla \phi_{r}|^2.
\endaligned\end{equation}
Under the assumption that \eqref{2Th-2} and \eqref{2Th-1}, we can derive from \eqref{3-Sec-18}
\begin{equation}\label{3-Sec-19}\aligned
0\leq&\Bigg[Q_g(M)-\frac{n-2}{4}\sqrt{\frac{2(n-2)}{n-1}}\Big(\int_M\Big|W\\
&+\frac{\sqrt{n}}{\sqrt{8}(n-2)}\mathring{{\rm Ric}} \mathbin{\bigcirc\mkern-15mu\wedge} g\Big|^{\frac{n}{2}}\Big)^{\frac{2}{n}}\Bigg]\Big(\int_M(|\mathring{R}_{ij}|\phi_{r})^{\frac{2n}{n-2}}\Big)^{\frac{n-2}{n}}\\
\leq&(1+\tilde{\epsilon}_3)\Big[\frac{1}{\tilde{\epsilon}_3}
+\frac{1}{\tilde{\epsilon}_2(1-\tilde{\epsilon}_2)}\Big]\int_M|\mathring{R}_{ij}|^2|\nabla \phi_{r}|^2\rightarrow 0
\endaligned\end{equation}
as $r\rightarrow \infty$, which shows that $M^n$ is Einstein.

{\it Case two:} When $4\leq n\leq6$ and $R\geq0$, for all $\epsilon_2,\epsilon_3$, we always have
\begin{equation}\label{3-Sec-20}
\frac{n-2}{4}-\frac{1+\epsilon_3}{1-\epsilon_2}<0.
\end{equation}
Therefore, under the condition \eqref{2Th-3}, there are $\overline{\epsilon}_2,\overline{\epsilon}_3$ small enough such that
\begin{equation}\label{3-Sec-22}\aligned
0\leq&\Bigg[Q_g(M)-\frac{1+\overline{\epsilon_3}}{1-\overline{\epsilon}_2}\sqrt{\frac{2(n-2)}{n-1}}\Big(\int_M\Big|W\\
&+\frac{\sqrt{n}}{\sqrt{8}(n-2)}\mathring{{\rm Ric}} \mathbin{\bigcirc\mkern-15mu\wedge} g\Big|^{\frac{n}{2}}\Big)^{\frac{2}{n}}\Bigg]\Big(\int_M(|\mathring{R}_{ij}|\phi_{r})^{\frac{2n}{n-2}}\Big)^{\frac{n-2}{n}}\\
\leq&\frac{1}{n-1}\Big[\frac{n-2}{4}-\frac{1+\overline{\epsilon}_3}{1-\overline{\epsilon}_2}\Big]
\int_MR|\mathring{R}_{ij}|^2\phi_{r}^2\\
&+(1+\overline{\epsilon}_3)\Big[\frac{1}{\overline{\epsilon}_3}
+\frac{1}{\overline{\epsilon}_2(1-\overline{\epsilon}_2)}\Big]\int_M|\mathring{R}_{ij}|^2|\nabla \phi_{r}|^2\rightarrow 0
\endaligned\end{equation}
as $r\rightarrow \infty$. Hence, $M^n$ is Einstein. In this case, \eqref{2Th-3} becomes
\begin{equation}\label{3-Sec-23}\aligned
\Big(\int_M|W|^{\frac{n}{2}}\Big)^{\frac{2}{n}}<\sqrt{\frac{n-1}{2(n-2)}}Q_g(M).
\endaligned\end{equation}
It is easy to check that for $n=4$ we have $\sqrt{\frac{n-1}{2(n-2)}}<\frac{\sqrt{6}}{4}$, for $n=5$ we have $\sqrt{\frac{n-1}{2(n-2)}}<\frac{2\sqrt{15}-4}{\sqrt{10}}$. This combining with Lemma \ref{2lem-4} shows that
$M^n(n=4,5)$ is of constant sectional curvature.

This completes the proof of Theorem \ref{thm1-2}.


\bibliographystyle{Plain}

\begin{thebibliography}{10}

\bibitem{Bach21}
R. Bach, Zur Weylschen Relativit\"{a}tstheorie und der Weylschen Erweiterung des Kr\"{u}mmungstensorbegriffs, Math. Z. 1921, 9: 110-135.

\bibitem{Catino2016} G. Catino,
Integral pinched shrinking Ricci solitons,
Adv. Math. 303 (2016), 279-294.

\bibitem{Chu2011}
Y. Chu, A rigidity theorem for complete noncompact Bach-flat manifolds, J. Geom. Phys. 61 (2011), 516-521.

\bibitem{Chu2012}
Y. Chu, P. Feng, Rigidity of complete noncompact Bach-flat $n$-manifolds, J. Geom. Phys. 62 (2012), 2227-2233.

\bibitem{Yuan2017}
Y. Fang, W. Yuan,
A sphere theorem for Bach-flat manifolds with positive constant scalar curvature,
arXiv:1704.06633

\bibitem{FuXiao2015}
H.-P, Fu, L.-Q, Xiao, Rigidity theorem for integral pinched
shrinking Ricci solitons, Monatsh. Math. 183 (2017), 487-494.

\bibitem{Huang2017}
G. Huang, Rigidity of Riemannian manifolds with positive scalar curvature, Ann. Glob. Anal. Geom. https://doi.org/10.1007/s10455-018-9600-x

\bibitem{Huang2017jmaa}
G. Huang, Integral pinched gradient shrinking $\rho$-Einstein solitons, J. Math. Anal. Appl. 451 (2017), 1045-1055.

\bibitem{HCL2018}
G. Huang, Y. Chen, X. Li,
Rigidity of Einstein metrics as critical points of some quadratic curvature functionals on complete manifolds, arXiv:1804.10748.

\bibitem{Kim2010}
S. Kim, Rigidity of noncompact complete Bach-flat manifolds, J. Geom. Phys. 60 (2010), 637-642.

\bibitem{Kim2011}
S. Kim, Rigidity of noncompact complete manifolds with harmonic curvature, Manuscripta Math. 135 (2011), 107-116.

\bibitem{MH2017}
B. Ma, G. Huang,
Rigidity of complete noncompact Riemannian manifolds with harmonic curvature, J. Geom. Phys. 124 (2018), 233--240.

\bibitem{MHLC2018}
B. Ma, G. Huang, X. Li, Y. Chen, Rigidity of Einstein metrics as critical points of quadratic curvature functionals on closed manifolds, Nonlinear Anal. 175 (2018), 237--248.

\bibitem{Schoen-yau1988}
R. Schoen, S.-T. Yau,  Conformally flat manifolds, Kleinian groups and scalar curvature, Invent. Math. 92 (1988), 47-71.



\end{thebibliography}

\end{document}